\documentclass[12pt]{amsart}
\usepackage{amscd,amssymb,graphicx,color,a4wide,cite, amsmath}
\usepackage{listings}
\usepackage{color}
\usepackage{booktabs}

\usepackage{hyperref}
\usepackage{url}

\definecolor{dkgreen}{rgb}{0,0.6,0}
\definecolor{gray}{rgb}{0.5,0.5,0.5}
\definecolor{mauve}{rgb}{0.58,0,0.82}

\usepackage{mathrsfs}

\usepackage{epstopdf}

\usepackage[all]{xy}
\newtheorem{theorem}{Theorem}
\newtheorem{lemma}[theorem]{Lemma}

\theoremstyle{definition}

\theoremstyle{remark}

\numberwithin{equation}{section}
\numberwithin{figure}{section}

\newcommand{\ZZ} {\mathbb{Z}}

\newcommand{\RR} {\mathbb{R}}

\newcommand{\PP} {\mathbb{P}}

\newcommand {\lra}  {\longrightarrow}

\newcommand {\Square}  {\operatorname{Sq}}

\def\mydate{\ifcase\month \or January\or February\or March\or
April\or May\or June\or July\or August\or September\or October\or 
November\or December\fi \space\number\day,\space\number\year}

\graphicspath{{images/}}

\newcommand\restr[2]{{
  \left.\kern-\nulldelimiterspace 
  #1 
  \vphantom{\big|} 
  \right|_{#2} 
  }}

\begin{document}

\title[Schoen's Calabi--Yau Threefolds]{Cohomology of the fixed point locus of an anti-symplectic involution on Schoen's Calabi--Yau}

\author{H\"ulya Arg\"uz}
\address{Labaratoure de Math\'ematiques de Versailles, Universit\'e de Versailles Saint-Quentin-en-Yvelines\\Versailles, $78000$\\France}
\email{nuromur-hulya.arguz@uvsq.fr}

\author{Thomas Prince}
\address{Mathematical Institute\\University of Oxford\\Woodstock Road\\Oxford\\OX$2$ $6$GG\\UK}
\email{thomas.prince@magd.ox.ac.uk}

\date{\today}

\begin{abstract}
We study the topology of a real Lagrangian in Schoen's Calabi--Yau threefold $X$ and compute its mod $2$ cohomology using two methods; first via a concrete Mayer--Vietoris calculation, then by an exact sequence relating the mod $2$ cohomology of the real Lagrangian to the cohomology of $X$. We conclude that these two methods agree. This in particular corrects a previous computation made by Casta\~{n}o-Bernard--Matessi.
\end{abstract}
\maketitle

\section{Introduction}
Schoen's Calabi--Yau threefold $X$ \cite{Sch88} is the fibred product of two rational elliptic surfaces over $\PP^1$, and is of particular interest from the point of view of the Strominger--Yau--Zaslow conjecture in mirror symmetry \cite{hosono1997mirror,GrossBB}. Kovalev has described a Lagrangian $3$-torus fibration on $X$ \cite{Kov03}, which has further been studied by Gross in \cite[\S$4$]{GrossBB} from the perspective of toric degenerations of Calabi--Yau toric complete intersections. In particular, Gross describes a singular torus fibration on $X$ using constructions given in \cite{GrossTopology}. 

Later on, Casta\~{n}o-Bernard--Matessi--Solomon studied anti-symplectic involutions on the total space of such torus fibrations in \cite{CBMS}. The fixed point set of such involutions are real Lagrangians. As these Lagrangians provide an algebro-geometric path to open Gromov-Witten invariants and the Fukaya category \cite{FOOO}, understanding their topology attracts significant interest from both symplectic and algebraic geometers. The mod $2$ cohomology of real Lagrangians is particularly widely studied in the literature due to its deep connections to real enumerative algebraic geometry ~\cite{bihan2003asymptotic,itenberg-topology,renaudineau2018bounding}. 

For the case of the real Lagrangian in Schoen's Calabi--Yau $X$, Casta\~{n}o-Bernard--Matessi outline the computation of the cohomology of a real Lagrangian in \cite[\S$4.2$]{CBM} via an argument based on the Mayer--Vietoris theorem. Here we first detail this Mayer-Vietoris computation by further investigating the topology of the torus fibration on Schoen's Calabi--Yau threefold, noting that our result differs slightly from the one obtained in \cite{CBM} (in which $H^1(\Sigma, \ZZ_2) \cong \ZZ_2^{34}$). We explain this discrepancy below, and give a detailed account of the Mayer--Vietoris calculation outlined in \cite{CBM}. Our main result is the following:
\begin{theorem}
\label{Thm main thm}
The mod $2$ cohomology groups of the real Lagrangian $\Sigma \subset X$ are given by
\[ H^0(\Sigma,\ZZ_2) = H^3(\Sigma,\ZZ_2) \cong \ZZ_2, ~~~ \mathrm{and} ~~~ H^1(\Sigma,\ZZ_2) = H^1(\Sigma,\ZZ_2) \cong \ZZ_2^{36}.  \]
\end{theorem}
We provide an alternative proof to this result using techniques developed in \cite{AP1}, where we give a general framework to compute the mod $2$ cohomology groups of real Lagrangians in a Calabi--Yau $X$, and its mirror $\breve{X}$, which admit (singular) torus fibrations
\[
f\colon X \to B, ~~~ \mathrm{and} ~~~ \breve{f}\colon \breve{X} \to B
\]
over a $\ZZ_2$ homology sphere $B$. We let $\iota$ denote the anti-symplectic involution on $X$ sending $x \mapsto -x$ in each fibre. Our main result in \cite{AP1} involves the comparison of the following maps:
\begin{enumerate}
	\item The map $\Square \colon D \mapsto D^2$, where $D \in H^1(B, R^1\breve{f}_\star\ZZ_2)$.
	\item The map $\beta \colon H^1(B, R^2f_\star\ZZ_2) \to H^2(B, R^1f_\star\ZZ_2)$; the connecting homomorphism introduced by Casta\~no-Bernard--Matessi in~\cite{CBM} associated to the short exact sequence 
\begin{equation}
\nonumber
0 \to R^1f_\star\ZZ_2\oplus\ZZ^2_2 \to \pi_\star\ZZ_2 \to R^2f_\star\ZZ_2 \to 0,
\end{equation}
where $\pi$ denotes the restriction of $f$ to the fixed locus $L_\RR$ of $\iota$.
\end{enumerate}
We prove in \cite[Theorem~$1.1$]{AP1} that the linear maps $\beta$ and $\Square$ coincide. Applying this result to the case of Schoen's Calabi--Yau, we verify Theorem \ref{Thm main thm}. This result also provides the first instance of a real Lagrangian in Calabi--Yau with Picard number greater than one which provides positive evidence for Conjecture~$1.3$ in \cite{APcoh}, stating that the cohomology of real Lagrangians in Calabi--Yau toric hypersurfaces are $2$-primary.

\subsection*{Acknowledgements}
TP was supported by a Fellowship by Examination at Magdalen College, Oxford and HA was supported by the Fondation Math\'{e}matiques Jacques Hadamard.

\section{A Real Lagrangian in Schoen's Calabi--Yau}
\label{sec:mayer_vietoris}
Take an elliptically fibred K$3$ surface over $\PP^1$, with 24 singular fibres
of Kodaira type $I_1$. Let $\gamma$ be a simple closed curve bounding a disk $D^2$ containing the image of $12$ singular fibers and not intersecting any other critical value on the base. Let $M^4$ be the preimage of $D^2$, and set $\overline{X} := M^4 \times T^2$, where $T^2$ denotes the standard $2$-torus. It is shown in \cite[\S$4$]{GrossBB} that Schoen's Calabi--Yau $X$ is diffeomorphic to a fibered coproduct of two copies of $\overline{X}$ along their boundary.

 There is a fibration on $\overline{X}$ over the solid torus $D^2 \times S^1$, obtained by taking the product of the fibration on $M^4$ with the standard $S^1$ fibration on $T^2$. This gives a torus fibration on $X$ over the three sphere $S^3$, expressed as the union of the solid tori 
 \[
 B_i=D^2 \times S^1,
 \]
 for $1\leq i \leq 2$, the bases of the fibrations on the two copies of $\overline{X}$ forming $X$. 

Consider an involution on $M^4$ obtained by restricting a fiber-preserving anti-symplectic, and anti-holomorphic involution on the K$3$ surface and take the involution on $T^2$ which preserves the fibres of the fibration on $S^1$. The product of these two involutions induces an anti-symplectic involution on $X$, and the fixed point locus $\Sigma$, which is naturally a real Lagrangian, is an $8$-to-$1$ cover
\[
\sigma \colon \Sigma \to B \cong S^3,
\] 
branched over $24$ circles, as described in \cite{CBM}. 
\subsection{Proof of Theorem~\ref{Thm main thm} via Mayer-Vietoris}
 We set $A_i := \sigma^{-1}(B_i)$ for each $i \in \{1,2\}$; we will apply the Mayer--Vietoris theorem to the decomposition $\Sigma = A_1 \cup A_2$.

First note that the fixed point locus on the fiber-preserving anti-symplectic involution on $M^4$ is obtained, as in \cite[\S$4.2$]{CBM}, as the dishoint uunion of a $2$-disc $E$ and a genus $4$ surface with three disks removed, denoted by $S$. The boundary components of $S$ are contained in the preimage of the boundary of $D^2$. Recalling that $B_i \cong D^2 \times S^1$ for each $i \in \{1,2\}$ we obtain that each $A_i$ is the disjoint union of two copies of $S\times S^1$, and two copies of $E \times S^1$.

We let $E^1_i \times S^1$, $S^1_i \times S^1$, $E^2_i \times S^1$, $S^2_i \times S^1$ denote the four components of $A_i$ for each $i \in \{1,2\}$. Following \cite{CBM, AP1} we may identify fibres of $\sigma$ with the half integral points in $\RR^2/\ZZ^2$ and $\RR^3/\ZZ^3$. Note that, similarly as in \cite[Corollary 1]{CBM}, since the monodromy always leaves one point on each fiber invariant, the fixed point locus $\Sigma$ has two connected components: one is homeomorphic to the base $S^3$ and the other is a $7$-to-$1$ branched cover over $S^3$.
For each $i \in \{1,2\}$, the boundary of $A_i$ consists of four two-dimensional tori (given by $\partial (E^j_i \times S^1)$ and $ \partial (S^j_i\times S^1)$ for $j \in \{1,2\}$ respectively). Following the constructions described in \cite[\S$4.1$]{CBM} and \cite[\S$4$]{GrossBB}, the components of $A_1\cap A_2$ can be indexed by $T_i$, for $0\leq i \leq 7$, where each $T_i$ is the $2$-torus defined as follows:
\begin{enumerate}
	\item $T_0$ bounds a copy of $E^1_i \times S^1$ in $A_i$ for each $i \in \{1,2\}$.
	\item $T_4$ bounds an the solid torus $E^2_1 \times S^1$ in $A_1$.
	\item $T_1$, $T_2$, and $T_3$ are boundary components of $S^1_1\times S^1 \subset A_1$.
	\item $T_5$, $T_6$, and $T_7$ are boundary components of $S^2_1 \times S^1 \subset A_1$.
	\item $T_3$ bounds the solid torus $E^2_2 \times S^1$ in $A_2$.
	\item $T_1$, $T_4$, and $T_5$ are boundary components of $S^1_2\times S^1 \subset A_2$.
	\item $T_2$, $T_6$, and $T_7$ are boundary components of $S^2_2\times S^1 \subset A_2$.
\end{enumerate}

Now consider the following part of the Mayer--Vietoris sequence for the decomposition $X = A_1 \cup A_2$:
\[
H_1(\Sigma,\ZZ_2) \to H_0(A_1\cap A_2,\ZZ_2) \to H_0(A_1,\ZZ_2)\oplus H_0(A_2,\ZZ_2) \to H_0(\Sigma,\ZZ_2) \to 0.
\]
Using the descriptions of $A_1$, $A_2$, and $X$ given above, this sequence has the form
\begin{equation}
\label{eq:explicit_1}
H_1(\Sigma,\ZZ_2) \to \ZZ_2^8 \to \ZZ^4_2 \oplus \ZZ^4_2 \to \ZZ_2^2 \to 0.
\end{equation}
The preceding part of the Mayer--Vietoris sequence has the form
\[
H_1(A_1\cap A_2,\ZZ_2) \to H_1(A_1,\ZZ_2)\oplus H_1(A_2,\ZZ_2) \to H_1(\Sigma,\ZZ_2).
\]
Computing these spaces, this sequence becomes
\begin{equation}
\label{eq:explicit_2}
\ZZ^{16}_2 \stackrel{\varphi}{\lra} \ZZ^{24}_2 \oplus \ZZ^{24}_2 \lra H_1(\Sigma,\ZZ_2).
\end{equation}
We let $\varphi = \varphi_1\oplus \varphi_2$ denote the components of the map $\varphi$ in \eqref{eq:explicit_2}. Note that it follows directly from \eqref{eq:explicit_1} and \eqref{eq:explicit_2} that $\dim H_1(\Sigma,\ZZ_2) = 34 + k$, where $k =\dim \ker(\varphi)$. We now describe a system of linear equations defining the kernel of $\varphi$.

Since $A_1 \cap A_2$ is the disjoint union of eight two-dimensional tori, the vector space $H_1(A_1\cap A_2,\ZZ_2)$ is generated by eight pairs $x_i,y_i \in H_1(T_i,\ZZ_2)$ for $i \in \{0,\ldots,7\}$.  Hence any element of $H_1(A_1\cap A_2,\ZZ_2)$ can be expressed as a sum of the form $\sum_{i=0}^7(a_ix_i+b_iy_i)$, where $a_i$ and $b_i$ are elements of $\ZZ_2$.

We choose the bases $\{x_i,y_i\}$ of each component of $H_1(A_1\cap A_2,\ZZ_2)$ such that $\varphi_1(y_i)$ generates the first homology group of the $S^1$ factor of the corresponding connected component of $A_1$ (recalling that this component is homeomorphic to either $E^j_1 \times S^1$ or $S^j_1 \times S^1$ for some $j \in \{1,2\}$) and $\varphi_2(x_i)$ generates the first homology group of the $S^1$ factor of the corresponding connected component of $A_2$. Moreover, the homology classes $\varphi_1(x_i)$ and $\varphi_2(y_i)$ are either meridian curves on a solid torus, or the classes of boundary circles of one of the three discs removed from the punctured genus four surface $S$.

\begin{lemma}
	\label{lem:equations}
	Consider a class $\sum_{i=0}^7(a_ix_i+b_iy_i)$, for $a_i$ and $b_i$ in $\ZZ_2$. This class is contained in the kernel of $\varphi$ if and only if the following equations hold:
\begin{enumerate}
	\item $a_0=b_0=0$.
	\item $b_4=0$.
	\item $b_1+b_2+b_3=0$ and $a_1=a_2=a_3$.
	\item $b_5+b_6+b_7=0$ and $a_5=a_6=a_7$.
	\item $a_3=0$.
	\item $a_1+a_4+a_5=0$ and $b_1=b_4=b_5$.
	\item $a_2+a_6+a_7=0$ and $b_2=b_6=b_7$.
\end{enumerate}
\end{lemma}
\begin{proof}
	These equations each follow directly from the above description of the $3$-manifolds bounded by each $T_i$ for $i \in \{0,\ldots, 7\}$. For example, considering the equations in item $(1)$, both $x_0$ and $y_0$ are the classes of longitudes of a solid torus and hence both $a_0$ and $b_0$ must vanish.
	
	Considering the equations in item $(3)$, we recall that $T_1$, $T_2$, and $T_3$ are the boundary components of a copy of $S \times S^1 \subset A_1$. The classes $y_1$, $y_2$, and $y_3$ are all mapped (by $\varphi_1$) to the multiples of the $S^1$ factor of $S^1_1 \times S^1$, and hence the component of $\varphi_1\big(\sum_{i=0}^7(a_ix_i+b_iy_i)\big)$ corresponding to this factor is equal to the class of $\varphi_1(b_1y_1+b_2y_2+b_3y_3) = b_1+b_2+b_3$, which must vanish. The classes $x_1$, $x_2$, and $x_3$ correspond to boundary components of the surface $S^1_1$. The only linear combination $\sum^3_{i=1}{a_ix_i}$ of these classes which are boundaries of $2$-chains in $S^1_1$ are multiples of $x_1+x_2+x_3$ (which is the class of the boundary of $S^1_1$ itself). That is, $\sum^3_{i=1}{a_ix_i} = \lambda (x_1+x_2+x_3)$ for some $\lambda$; hence $a_1=a_2=a_3$. Similar computations verify the equations in items $(4)$, $(6)$, and $(7)$.
\end{proof}

Taken as equations over the integers, the equations described in Lemma~\ref{lem:equations} have no non-zero solutions, and hence $H^1(\Sigma,\ZZ) \cong \ZZ^{34}$ (in agreement with the calculation of \cite[\S$4.2$]{CBM}). However, over $\ZZ_2$, we find a $2$-dimensional space of solutions generated by:
\begin{enumerate}
	\item $a_i = 0$ for all $i \in \{0,\ldots,7\}$, $b_0=b_1=b_4=b_5=0$, and $b_2=b_3=b_6=b_7=1$,
	\item $b_i = 0$ for all $i \in \{0,\ldots,7\}$, $a_0=a_1=a_2=a_3=0$, and $a_4=a_5=a_6=a_7=1$.
\end{enumerate}

That is, we find that $k=2$, $H_1(\Sigma,\ZZ_2) \cong \ZZ_2^{36}$, and hence $H^1(\Sigma,\ZZ_2) \cong \ZZ_2^{36}$.

\subsection{Proof of Theorem~\ref{Thm main thm} via the squaring map}
Consider the following exact sequence of sheaf cohomology groups introduced in \cite{CBM}:
\begin{equation}
\label{eq:les_section}
0 \to H^1(B,R^1f_\star\ZZ_2) \to H^1(B,\sigma_\star\ZZ_2) \to H^1(B,R^2f_\star\ZZ_2) \overset{\beta}{\to} H^2(B,R^1f_\star\ZZ_2).
\end{equation}
Noting that Schoen's Calabi--Yau $X$ can be described as a complete intersection in a toric variety \cite{Sch88}, it is simply-connected by the Lefschetz hyperplane theorem; in this case, $H^1(B,R^1\breve{f}_\star\ZZ_2) \cong H_2
(X, \ZZ_2)$. Moreover, the mirror to Schoen's Calabi--Yau $X$ is itself homeomorphic to $X$ (see \cite[\S$4$]{GrossBB}). Hence the ranks of $H^1(B,R^1\breve{f}_\star\ZZ_2)$ and $H^1(B,R^2f_\star\ZZ_2)$ agree and \cite[Theorem 1.1]{AP1} implies that the rank of the connecting homomorphism
\begin{equation}
\label{eq:beta}
\beta \colon H^1(B,R^2f_\star\ZZ_2) \to H^2(B,R^1f_\star\ZZ_2)
\end{equation}
is equal to the rank of the squaring map $\Square \colon D \mapsto D^2$ from $H^2(X,\ZZ_2)$ to $H^4(X,\ZZ_2)$ which takes an element $D \in H^2(X,\ZZ_2)$ to $D\smile D \in H^4(X,\ZZ_2)$. Note that this map is linear over $\ZZ_2$.

As discussed in the introduction, Schoen's Calabi--Yau is the fibre product of a pair of rational elliptic surfaces $E_1$ and $E_2$
	\begin{equation}
	\label{eq:square}
	\xymatrix{
	X \ar^{\pi_2}[r] \ar^{\pi_1}[d] & E_2 \ar[d] \\
	E_1 \ar[r] & \PP^1.
	}
	\end{equation}
	In particular, $H^2(X,\ZZ_2) \cong \ZZ_2^{19}$ is spanned by vectors in the $10$ dimensional subspaces $\pi_1^{\star}(H^2(E_1,\ZZ_2)) \cong \ZZ_2^{10}$ and $\pi_2^{\star}(H^2(E_2,\ZZ_2)) \cong \ZZ_2^{10}$. These spaces intersect in the class of a fibre of either composition $X \to \PP^1$ in \eqref{eq:square}.
	
	Since $(\pi_i^\star(\alpha))^2 = \pi_i^\star(\alpha^2)$, for any class $\alpha \in H^2(E_i,\ZZ_2)$ and $i \in \{1,2\}$, the image of $\Square$ is spanned by $\pi_1^\star(H^4(E_1,\ZZ_2)) \cong \ZZ_2$ and $\pi_2^\star(H^4(E_2,\ZZ_2)) \cong \ZZ_2$. That is, the rank of $\Square$ is equal to $2$ and, by \cite[Theorem~$1.1$]{AP1}, the rank of $\beta$ is also equal to $2$.
In particular, we can conclude from \eqref{eq:les_section} that $H^1(B,\sigma_\star\ZZ_2)$ has dimension $19+(19-2)=36$. Hence, from the Leray spectral sequence for $\sigma$, we have that $H^1(B,\sigma_\star\ZZ_2) \cong H^1(\Sigma, \ZZ_2) \cong \ZZ_2^{36}$, verifying Theorem~\ref{Thm main thm}.

\bibliographystyle{plain}
\bibliography{bibliography}

\end{document}